\documentclass[12pt]{article}
\usepackage{amsmath,amsfonts,amssymb,amsthm}
\usepackage{color}
\usepackage{colortbl}
\usepackage{array}
\usepackage{mathrsfs}
\usepackage{dcolumn}
\usepackage{longtable}
\usepackage{hhline}
\usepackage{bm}

\newcommand{\braced}[2]{\genfrac{\{}{\}}{0pt}{0}{#1}{#2}}

\theoremstyle{plain}
\newtheorem{thm}{Theorem}[section]

\newtheorem{cor}[thm]{Corollary}

\theoremstyle{definition}
\newtheorem{defn}[thm]{Definition}

\title{\bf On Generalized Multi Poly-Euler and\\ Multi Poly-Bernoulli Polynomials}
\author{
{\large Roberto B. Corcino, Hassan Jolany, Cristina B. Corcino and Takao Komatsu
}
}

\date{}

\begin{document}

\maketitle

\begin{abstract}
In this paper, we establish more identities of generalized multi poly-Euler polynomials with three parameters and obtain a kind of symmetrized generalization of the polynomials. Moreover, generalized multi poly-Bernoulli polynomials are defined using multiple polylogarithm and derive some properties parallel to those of poly-Bernoulli polynomials. These are generalized further using the concept of Hurwitz-Lerch multiple zeta values.    


\end{abstract}

\section{Introduction}
The Euler numbers, denoted by $E_n$, are usually introduced as the coefficients of the generating function
$$\frac{1}{\cosh t}=\frac{2e^t}{e^{2t}+1}=\sum_{n=0}^{\infty}E_n\frac{t^n}{n!}.$$
Eventually, these numbers have been generalized in polynomial form as
$$\frac{2e^{xt}}{e^t+1}=\sum_{n=0}^{\infty}E_n(x)\frac{t^n}{n!}$$
where $E_n(x)$ denote the so-called {\it Euler polynomials}. 

\smallskip
Euler numbers and polynomials have a rich literature in the history of mathematics where bunch of identities and properties have been established including their mathematical and physical applications. These numbers and polynomials have close connections with Bernoulli numbers and polynomials, particularly, in the structures of their properties and generalizations. In fact, in almost every property and generalization of Bernoulli numbers and polynomials there corresponds property and generalization for Euler numbers and polynomials. For instance, Kaneko \cite{Kaneko} introduced the poly-Bernoulli numbers $B^{(k)}_n$ by means of the following exponential generating function
$$\frac{{\rm Li}_k(1-e^{-x})}{1-e^{-x}}=\sum_{n=0}^{\infty}B^{(k)}_n\frac{x^n}{n!}$$  
where
$${\rm Li}_k(z)=\sum_{n=1}^{\infty}\frac{z^n}{n^k},$$
while Ohno and Sasaki \cite{Sasaki} defined poly-Euler numbers as
$$\frac{{\rm Li}_k(1-e^{-4t})}{4t\cosh t}=\sum_{n=0}^{\infty}E^{(k)}_n\frac{t^n}{n!}$$
which have been recently extended by H. Jolany et al. \cite{Jolany3} in polynomial form as
\begin{equation}\label{polyform}
\frac{2{\rm Li}_k(1-e^{-t})}{1+e^{t}}e^{xt}=\sum_{n=0}^{\infty}E^{(k)}_n(x)\frac{t^n}{n!}.
\end{equation}

\smallskip
It is worth-mentioning that the above generalization of Kaneko has been generalized further by Cenkci and Young \cite{Cenkci} using the concept of Hurwitz-Lerch zeta function $\Phi (z,s,a)$ as follows
\begin{equation}\label{HpolyB}
\Phi (1-e^{-t},k,a)=\sum_{n=0}^{\infty}B_{n,a}^{(k)}\frac{t^n}{n!}
\end{equation}
where
\begin{equation}\label{HLerch}
\Phi (z,k,a)=\sum_{n=0}^{\infty}\frac{z^n}{(n+a)^k}.
\end{equation}
The numbers $B_{n,a}^{(k)}$ are called the Hurwitz type poly Bernoulli numbers. These numbers have been shown to have explicit formula
\begin{equation}\label{HLExplicit}
B_{n,a}^{(k)}=(-1)^n\sum_{m=0}^{n}\frac{(-1)^mm!S(n,m)}{(m+a)^k}
\end{equation}
where $S(n,m)$ denotes the Stirling numbers of the second kind. A parallel version of generalization for Euler numbers is still to be done. However, it is expected that the structure is quite complicated. For instance, one may define it as
\begin{equation}\label{HLEuler}
\frac{2(1-e^{-t})\Phi (1-e^{-t}, k, a)}{1+e^t}=\sum_{n=0}^{\infty}E_{n,a}^{(k)}\frac{t^n}{n!}.
\end{equation}
The numbers $E_{n,a}^{(k)}$ may be called the Hurwitz type poly Euler numbers. One can easily derive the explicit formula for $E_{n,a}^{(k)}$ as follows
\begin{equation}\label{HLExplicit1}
E_{j+1,a}^{(k)}=\sum_{n=j}^{\infty}\sum_{r=j}^{n}\sum_{q=j}^r\sum_{m=0}^{n-r}\frac{(-1)^{n-r+m+q+1}\binom{j+1}{q-n+j+1,\;r-q,\;n-r}(n-j)^{r-q}m!S(n-r,m)}{(m+a)^k}.
\end{equation}
We notice that the explicit formula of $E_{n,a}^{(k)}$ is more complicated than that of $B_{n,a}^{(k)}$. 

\smallskip
Another form of generalization of Bernoulli numbers has been defined by Imatomi et al. \cite{Imatomi} in terms of multiple polylogarithm as follows
\begin{equation}\label{MultiBernoulli}
\frac{{\rm Li}_{(k_1, k_2,\ldots, k_r)}(1-e^{-t})}{1-e^{-t}}=\sum_{n=0}^{\infty}{B}^{(k_1, k_2,\ldots, k_r)}_n\frac{t^n}{n!}
\end{equation}
where 
\begin{equation}\label{multipolylog}
{\rm Li}_{(k_1,k_2,\ldots, k_r)}(z)=\sum_{ 0< m_1< m_2<\ldots < m_r }\frac{z^{m_r}}{m_1^{k_1} m_2^{k_2}\ldots m_r^{k_r}}.
\end{equation}
These numbers possess respectively the following recurrence relation and explicit formula
\begin{align}
{B}^{(k_1, k_2,\ldots, k_r)}_n&=\frac{1}{n+1}\left({B}^{(k_1-1, k_2,\ldots, k_r)}_n-\sum_{m=1}^{n-1}\binom{n}{m-1}{B}^{(k_1, k_2,\ldots, k_r)}_m\right)\\
{B}^{(k_1, k_2,\ldots, k_r)}_n&=(-1)^n\sum_{n+1\geq m_1>m_2>\ldots >m_r>0}\frac{(-1)^{m_1-1}(m_1-1)!S(n,m_1-1)}{m_1^{k_1}m_1^{k_2}\ldots m_r^{k_r}}.
\end{align}
Parallel to the above generalization is the generalized multi poly-Euler polynomials which are denoted by ${E}^{(k_1, k_2,\ldots, k_r)}_n(x;a,b,c)$. These polynomials have been introduced in \cite{CorJol1} by means of the above multiple poly-logarithm, also known as multiple zeta values. More precisely, we have
\begin{equation}\label{multipolyeuler}
\frac{2{\rm Li}_{(k_1, k_2,\ldots, k_r)}(1-(ab)^{-t})}{(a^{-t}+b^{t})^r}c^{rxt}=\sum_{n=0}^{\infty}{E}^{(k_1, k_2,\ldots, k_r)}_n(x;a,b,c)\frac{t^n}{n!}.
\end{equation}
When $r=1$, (\ref{multipolyeuler}) boils down to the generalized poly-Euler polynomials with three parameters $a, b, c$. Moreover, when $c=e$, (\ref{multipolyeuler}) reduces to the multi poly-Euler polynomials with two parameters $a, b$. These special cases have been discussed intensively in \cite{CorJol1, CorJol2}.

\smallskip
Some properties of generalized multi poly-Euler polynomials ${E}^{(k_1, k_2,\ldots, k_r)}_n(x;a,b,c)$ are established in \cite{CorJol1} which include the following identities
\begin{align}
{E}^{(k_1, k_2, \ldots, k_r)}_n(x;a,b,c)&=\sum_{i=0}^n\binom{n}{i}(r\ln c)^{n-i}{E}^{(k_1, k_2, \ldots, k_r)}_i(a,b)x^{n-i}\label{prop1}\\
{E}^{(k_1, k_2, \ldots, k_r)}_n(x;a,b,c)&=(\ln a+\ln b)^n{E}^{(k_1, k_2, \ldots, k_r)}_n\left(\frac{rx\ln c+\ln a}{\ln a+\ln b}\right)\label{prop2}
\end{align}
\begin{align*}
\frac{d}{dx}{E}^{(k_1, k_2, \ldots, k_r)}_{n+1}(x;a,b,c)&=(n+1)(r\ln c){E}^{(k_1, k_2, \ldots, k_r)}_{n}(x;a,b,c)\\
{E}^{(k_1, k_2, \ldots, k_r)}_{n}(x+y;a,b,c)&=\sum_{i=0}^n\binom{n}{i}(r\ln c)^{n-i}{E}^{(k_1, k_2, \ldots, k_r)}_{i}(x;a,b,c)y^{n-i}.\label{prop4}
\end{align*}
Furthermore, an explicit formula for ${E}^{(k_1, k_2,\ldots, k_r)}_n(x;a,b,c)$ is given by
\begin{equation}\label{explicit1}
{E}^{(k_1,k_2,\ldots, k_r)}_n(x;a,b,c)=\sum_{i=0}^n\sum_{ 0< m_1< m_2<\ldots < m_r \atop c_1+c_2+\ldots=r}\mathcal{J}(m_1,m_2, \ldots, m_r)
\end{equation}
where 
$$\mathcal{J}(m_1,m_2, \ldots, m_r)=\sum_{j=0}^{m_r}\frac{2(rx\ln c-j\ln ab)^{n-i}r!(-1)^{j+s}(s\ln ab+r\ln a)^i\binom{m_r}{j}\binom{n}{i}}{(c_1!c_2!\ldots)(m_1^{k_1} m_2^{k_2}\ldots m_r^{k_r})}$$
with $s=c_1+2c_2+\ldots$.
When $r=1$, these identities reduce to
\begin{align*}
{E}^{(k)}_n(x;a,b,c)&=\sum_{i=0}^n\binom{n}{i}(\ln c)^{n-i}{E}^{(k)}_i(a,b)x^{n-i}\\
{E}^{(k)}_n(x;a,b,c)&=(\ln a+\ln b)^n{E}^{(k)}_n\left(\frac{x\ln c+\ln a}{\ln a+\ln b}\right)\\
\frac{d}{dx}{E}^{(k)}_{n+1}(x;a,b,c)&=(n+1)(\ln c){E}^{(k)}_{n}(x;a,b,c)\\
{E}^{(k)}_{n}(x+y;a,b,c)&=\sum_{i=0}^n\binom{n}{i}(\ln c)^{n-i}{E}^{(k)}_i(x;a,b,c)y^{n-i}
\end{align*}
which are properties of generalized poly-Euler polynomials (see \cite{CorJol1}).

\smallskip
In this paper, some identities of ${E}^{(k_1, k_2,\ldots, k_r)}_n(x;a,b,c)$ related to Stirling numbers of the second kind are established and certain symmetrized generalization for ${E}^{(k_1, k_2,\ldots, k_r)}_n(x;a,b,c)$ is obtained. On the other hand, generalized multi poly-Bernoulli polynomials are defined and some properties of these polynomials are established parallel to those of the poly-Bernoulli polynomials. Moreover, certain generalization of multi poly-Bernoulli numbers is defined in terms of generalized Hurwitz-Lerch multiple zeta values.

\section{Generalized Multi Poly-Euler Polynomials and Stirling Numbers}
The generalized poly-Euler polynomials with three parameters $a$, $b$ and $c$ are defined in \cite{CorJol1} as follows
\begin{equation}\label{genpolyeuler1}
\frac{2{\rm Li}_k(1-(ab)^{-t})}{a^{-t}+b^{t}}c^{xt}=\sum_{n=0}^{\infty}{E}^{(k)}_n(x;a,b,c)\frac{t^n}{n!}.
\end{equation}
Some identities on generalized poly-Euler polynomials are expressed in terms of Stirling numbers of the second kind. Such identities have appeared in Theorem 2.6 of \cite{CorJol1} but with $c=e$. More precisely, we have the following theorem.
\begin{thm}\label{add1} {\rm \cite{CorJol1}} 
The generalized poly-Euler polynomials ${E}^{(k)}_n(x;a,b)$ satisfy the following explicit formulas 
\begin{align}
{E}^{(k)}_{n}(x;a,b)&=\sum_{m=0}^{\infty}\sum_{l=m}^n\braced{l}{m}\binom{n}{l}{E}^{(k)}_{n-l}(-m;a,b)(x)^{(m)}\label{add1eq1}\\
{E}^{(k)}_{n}(x;a,b)&=\sum_{m=0}^{\infty}\sum_{l=m}^n\braced{l}{m}\binom{n}{l}{E}^{(k)}_{n-l}(0;a,b)(x)_{m}\label{add1eq2}\\
{E}^{(k)}_{n}(x;a,b)&=\sum_{m=0}^{\infty}\binom{n}{m}\sum_{l=0}^{n-m}\frac{\binom{n-m}{l}}{\binom{l+s}{l}}\braced{l+s}{s}{E}^{(k)}_{n-m-l}(0;a,b)B^{(s)}_m(x)\label{add1eq3}\\
{E}^{(k)}_{n}(x;a,b)&=\sum_{m=0}^{\infty}\frac{\binom{n}{m}}{(1-\lambda)^s}\sum_{j=0}^s\binom{s}{j}(-\lambda)^{s-j}{E}^{(k)}_{n-m}(j;a,b)H^{(s)}_m(x;\lambda),\label{add1eq4}
\end{align}
where $(x)^{(n)}=x(x+1)\ldots (x+n-1)$, $(x)_n=x(x-1)\ldots (x-n+1)$,
$$\left(\frac{t}{e^t-1}\right)^se^{xt}=\sum_{n=0}^{\infty}B^{(s)}_n(x)\frac{t^n}{n!}\mbox{  and  }\left(\frac{1-\lambda}{e^t-\lambda}\right)^se^{xt}=\sum_{n=0}^{\infty}H^{(s)}_n(x;\lambda)\frac{t^n}{n!}.$$
\end{thm}

Here, we derive some identities for ${E}^{(k_1, k_2,\ldots, k_r)}_n(x;a,b,c)$ which are parallel to those in Theorem \ref{add1} . The first such identity is given in the following theorem.
\begin{thm}\label{add2}
\begin{equation}\label{add1eq11}
{E}^{(k_1, k_2,\ldots, k_r)}_n(x;a,b,c)=\sum_{m=0}^{\infty}\sum_{l=m}^n(r\log c)^l\braced{l}{m}\binom{n}{l}{E}^{(k_1, k_2,\ldots, k_r)}_{n-l}(-m\log c;a,b)(x)^{(m)}.
\end{equation}
\begin{proof}
Note that (\ref{multipolyeuler}) can be written as
$$\sum_{n=0}^{\infty}{E}^{(k_1, k_2,\ldots, k_r)}_n(x;a,b,c)\frac{t^n}{n!}=\frac{2{\rm Li}_{(k_1, k_2,\ldots, k_r)}(1-(ab)^{-t})}{(a^{-t}+b^{t})^r}(1-(1-e^{-rt\log c}))^{-x}$$
Using Newton's Binomial Theorem, we have
$$\sum_{n=0}^{\infty}{E}^{(k_1, k_2,\ldots, k_r)}_n(x;a,b)\frac{t^n}{n!}=\frac{2{\rm Li}_{(k_1, k_2,\ldots, k_r)}(1-(ab)^{-t})}{(a^{-t}+b^{t})^r}\sum_{m=0}^{\infty}\binom{x+m-1}{m}(1-e^{-rt\log c})^m$$
\begin{align*}
&=\sum_{m=0}^{\infty}(x)^{(m)}\frac{(e^{rt\log c}-1)^m}{m!}\frac{2{\rm Li}_{(k_1, k_2,\ldots, k_r)}(1-(ab)^{-t})}{(a^{-t}+b^{t})^r}e^{-mrt\log c}\\
&=\sum_{m=0}^{\infty}(x)^{(m)}\left(\sum_{n=0}^{\infty}\braced{n}{m}\frac{(rt\log c)^n}{n!}\right)\left(\sum_{n=0}^{\infty}{E}^{(k_1, k_2,\ldots, k_r)}_n(-mr\log c;a,b)\frac{t^n}{n!}\right)\\
&=\sum_{n=0}^{\infty}\left\{\sum_{m=0}^{\infty}\sum_{l=m}^{n}(r\log c)^l\braced{l}{m}\binom{n}{l}{E}^{(k_1, k_2,\ldots, k_r)}_{n-l}(-mr\log c;a,b)(x)^{(m)}\right\}\frac{t^n}{n!}
\end{align*}
Comparing coefficients completes the proof of (\ref{add1eq11}).
\end{proof}
\end{thm}

In particular, when $c=e$ and $r=1$, (\ref{add1eq11}) yields (\ref{add1eq1}). For the generalization of (\ref{add1eq2}), we have the following theorem.

\begin{thm}\label{add3}
\begin{equation}\label{add1eq22}
{E}^{(k_1, k_2,\ldots, k_r)}_n(x;a,b,c)=\sum_{m=0}^{\infty}\sum_{l=m}^n(r\log c)^l\braced{l}{m}\binom{n}{l}{E}^{(k_1, k_2,\ldots, k_r)}_{n-l}(0;a,b)(x)^{(m)}.
\end{equation}
\begin{proof}
Note that (\ref{multipolyeuler}) can be written as
$$\sum_{n=0}^{\infty}{E}^{(k_1, k_2,\ldots, k_r)}_n(x;a,b,c)\frac{t^n}{n!}=\frac{2{\rm Li}_{(k_1, k_2,\ldots, k_r)}(1-(ab)^{-t})}{(a^{-t}+b^{t})^r}((e^{rt\log c}-1)+1)^{x}$$
Using Newton's Binomial Theorem, we have
$$\sum_{n=0}^{\infty}{E}^{(k_1, k_2,\ldots, k_r)}_n(x;a,b)\frac{t^n}{n!}=\frac{2{\rm Li}_{(k_1, k_2,\ldots, k_r)}(1-(ab)^{-t})}{(a^{-t}+b^{t})^r}\sum_{m=0}^{\infty}\binom{x}{m}(e^{rt\log c}-1)^m$$
\begin{align*}
&=\sum_{m=0}^{\infty}(x)_m\frac{(e^{rt\log c}-1)^m}{m!}\frac{2{\rm Li}_{(k_1, k_2,\ldots, k_r)}(1-(ab)^{-t})}{(a^{-t}+b^{t})^r}\\
&=\sum_{m=0}^{\infty}(x)_m\left(\sum_{n=0}^{\infty}\braced{n}{m}\frac{(rt\log c)^n}{n!}\right)\left(\sum_{n=0}^{\infty}{E}^{(k_1, k_2,\ldots, k_r)}_n(0;a,b)\frac{t^n}{n!}\right)\\
&=\sum_{n=0}^{\infty}\left\{\sum_{m=0}^{\infty}\sum_{l=m}^{n}(r\log c)^l\braced{l}{m}\binom{n}{l}{E}^{(k_1, k_2,\ldots, k_r)}_{n-l}(0;a,b)(x)_m\right\}\frac{t^n}{n!}
\end{align*}
Comparing coefficients completes the proof of (\ref{add1eq22}).
\end{proof}
\end{thm}

\begin{thm}\label{add4}
\begin{equation}\label{add1eq33}
{E}^{(k_1, k_2,\ldots, k_r)}_n(x;a,b,c)=\sum_{m=0}^{\infty}\binom{n}{m}\sum_{l=0}^{n-m}\frac{\binom{n-m}{l}}{\binom{l+s}{l}}\braced{l+s}{s}{E}^{(k_1, k_2,\ldots, k_r)}_{n-m-l}(0;a,b)B^{(s)}_m(xr\log c).
\end{equation}
\begin{proof}
Note that (\ref{multipolyeuler}) can be written as
$$\sum_{n=0}^{\infty}{E}^{(k_1, k_2,\ldots, k_r)}_n(x;a,b,c)\frac{t^n}{n!}=\left(\frac{(e^t-1)^s}{s!}\right)\left(\frac{t^se^{xrt\log c}}{(e^t-1)^s}\right)\left(\frac{2{\rm Li}_{(k_1, k_2,\ldots, k_r)}(1-(ab)^{-t})}{(a^{-t}+b^{t})^r}\right)\frac{s!}{t^s}$$
\begin{align*}
&=\left(\sum_{n=0}^{\infty}\braced{n+s}{s}\frac{t^{n+s}}{(n+s)!}\right)\left(\sum_{m=0}^{\infty}B^{(s)}_m(xr\log c)\frac{t^m}{m!}\right)\left(\sum_{n=0}^{\infty}{E}^{(k_1, k_2,\ldots, k_r)}_n(0;a,b)\frac{t^n}{n!}\right)\frac{s!}{t^s}\\
&=\left(\sum_{n=0}^{\infty}\braced{n+s}{s}\frac{t^{n+s}}{(n+s)!}\right)\left(\sum_{m=0}^{\infty}\sum_{n=m}^{\infty}B^{(s)}_m(xr\log c)\frac{t^m}{m!}{E}^{(k_1, k_2,\ldots, k_r)}_{n-m}(0;a,b)\frac{t^{n-m}}{(n-m)!}\right)\frac{s!}{t^s}\\
&=\sum_{m=0}^{\infty}\left\{\sum_{n=m}^{\infty}\sum_{l=0}^{n-m}\braced{l+s}{s}\frac{t^{l+s}}{(l+s)!}B^{(s)}_m(xr\log c){E}^{(k_1, k_2,\ldots, k_r)}_{n-m-l}(0;a,b)\frac{t^{n-m-l}}{(n-m-l)!}\frac{t^m}{m!}\frac{s!}{t^s}\right\}\\
&=\sum_{n=0}^{\infty}\left\{\sum_{m=0}^{n}\binom{n}{m}\sum_{l=0}^{n-m}\frac{\binom{n-m}{l}}{\binom{l+s}{s!}}\braced{l+s}{s}{E}^{(k_1, k_2,\ldots, k_r)}_{n-m-l}(0;a,b)B^{(s)}_m(xr\log c)\right\}\frac{t^n}{n!}
\end{align*}
Comparing coefficients completes the proof of (\ref{add1eq33}).
\end{proof}
\end{thm}

\begin{thm}\label{add5}
\begin{equation}\label{add1eq44}
{E}^{(k_1, k_2,\ldots, k_r)}_n(x;a,b,c)=\sum_{m=0}^{n}\frac{\binom{n}{m}}{(1-\lambda)^s}\sum_{j=0}^s\binom{s}{j}(-\lambda)^{s-j}{E}^{(k_1, k_2,\ldots, k_r)}_{n-m}(j;a,b)H^{(s)}_m(xr\log c;\lambda).
\end{equation}
\begin{proof}
Note that (\ref{multipolyeuler}) can be written as
$$\sum_{n=0}^{\infty}{E}^{(k_1, k_2,\ldots, k_r)}_n(x;a,b,c)\frac{t^n}{n!}=\left(\frac{(1-\lambda)^s}{(e^t-\lambda)^s}e^{xrt\log c}\right)\left(\frac{(e^t-\lambda)^s}{(1-\lambda)^s}\right)\left(\frac{2{\rm Li}_{(k_1, k_2,\ldots, k_r)}(1-(ab)^{-t})}{(a^{-t}+b^{t})^r}\right)$$
\begin{align*}
&=\left(\sum_{n=0}^{\infty}H^{(s)}_n(xr\log c;\lambda)\frac{t^n}{n!}\right)\left(\sum_{j=0}^{s}\binom{s}{j}(-\lambda)^{s-j}\frac{2{\rm Li}_{(k_1, k_2,\ldots, k_r)}(1-(ab)^{-t})}{(a^{-t}+b^{t})^r}e^{jt}\right)\\
&=\frac{1}{(1-\lambda)^s}\sum_{j=0}^{s}\binom{s}{j}(-\lambda)^{s-j}\left(\sum_{n=0}^{\infty}H^{(s)}_n(xr\log c;\lambda)\frac{t^n}{n!}\right)\left(\sum_{n=0}^{\infty}{E}^{(k_1, k_2,\ldots, k_r)}_{n}(j;a,b)\frac{t^{n}}{n!}\right)\\
&=\frac{1}{(1-\lambda)^s}\sum_{j=0}^{s}\binom{s}{j}(-\lambda)^{s-j}\sum_{n=0}^{\infty}\sum_{m=0}^{n}\binom{n}{m}H^{(s)}_m(xr\log c;\lambda){E}^{(k_1, k_2,\ldots, k_r)}_{n-m}(j;a,b)\frac{t^n}{n!}\\
&=\sum_{n=0}^{\infty}\left(\sum_{m=0}^{n}\frac{\binom{n}{m}}{(1-\lambda)^s}\sum_{j=0}^{s}\binom{s}{j}(-\lambda)^{s-j}{E}^{(k_1, k_2,\ldots, k_r)}_{n-m}(j;a,b)H^{(s)}_m(xr\log c;\lambda)\right)\frac{t^n}{n!}.
\end{align*}
Comparing coefficients completes the proof of (\ref{add1eq44}).
\end{proof}
\end{thm}

\smallskip
The next theorem contains an identity which is obtained by making use of the following differential formula for the generalized poly-logarithm (Hamahata and Masubuchi, Integers) 
\begin{equation}\label{df}
\frac{d}{d z}{\rm Li}_{k_1,\dots,k_r}(z)=
\begin{cases}  
\frac{1}{z}{\rm Li}_{k_1,\dots,k_{r-1},k_r-1}(z)&\text{if $k_r>1$};\\
\frac{1}{1-z}{\rm Li}_{k_1,\dots,k_{r-1}}(z)&\text{if $k_r=1$}. 
\end{cases} 
\end{equation}

\smallskip
\begin{thm}\label{add6}
If $k_r>1$, then 
\begin{align*} 
E_{n+1}^{(k_1,\dots,k_r)}(x)&-r\log(a c^x)E_n^{(k_1,\dots,k_r)}(x)\\
&=r\sum_{\nu=0}^\infty\binom{n}{\nu}\sum_{j=0}^\infty\bigl(-j\log(a b)\bigr)^{n-\nu}E_\nu^{(k_1,\dots,k_r)}(x)\\
&-\log(a b)\sum_{\nu=0}^\infty\binom{n}{\nu}\sum_{j=0}^\infty\bigl(j\log(a b)\bigr)^{n-\nu}E_\nu^{(k_1,\dots,k_{r-1},k_r-1)}(x)\,. \qquad\qquad\quad\;\;\;
\end{align*} 
If $k_r=1$, then 
\begin{align*} 
E_{n+1}^{(k_1,\dots,k_r)}(x)&-r\log(a c^x)E_n^{(k_1,\dots,k_r)}(x)\\
&=r\sum_{\nu=0}^\infty\binom{n}{\nu}\sum_{j=0}^\infty\bigl(-j\log(a b)\bigr)^{n-\nu}E_\nu^{(k_1,\dots,k_r)}(x)
+\log(a b)E_n^{(k_1,\dots,k_{r-1})}(x)\,. 
\end{align*} 
\begin{proof}
We differentiate both sides of 
$$
2(a c^x)^{r t}{\rm Li}_{k_1,\dots,k_r}\bigl(1-(a b)^{-t}\bigr)=\bigl(1+(a b)^t\bigr)^r\sum_{n=0}^\infty E_n^{(k_1,\dots,k_r)}(x)\frac{t^n}{n!}
$$ 
with respect to $t$ using (\ref{df}).  If $k_1\ne 1$, then 
\begin{align*} 
2 r(a c^x)^{r t}\log(a c^x)&{\rm Li}_{k_1,\dots,k_r}\bigl(1-(a b)^{-t}\bigr)\\
+2(a c^x)^{r t}&\frac{-(a b)^{-t}\log(a b)}{1-(a b)^{-t}}{\rm Li}_{k_1,\dots,k_{r-1},k_r-1}\bigl(1-(a b)^{-t}\bigr)\\
&=r\bigl(1+(a b)^t\bigr)^{r-1}(a b)^t\log(a b)\sum_{n=0}^\infty E_n^{(k_1,\dots,k_r)}(x)\frac{t^n}{n!}\\
&\;\;\;\;+\bigl(1+(a b)^t\bigr)^r\sum_{n=1}^\infty E_n^{(k_1,\dots,k_r)}(x)\frac{t^{n-1}}{(n-1)!}\,. 
\end{align*} 
Dividing both sides by $\bigl(1+(a b)^t\bigr)^r$, we have 
\begin{align*} 
r\log(a c^x)\frac{2(a c^x)^{r t}}{\bigl(1+(a b)^t\bigr)^r}&{\rm Li}_{k_1,\dots,k_r}\bigl(1-(a b)^{-t}\bigr)\\
+\frac{\log(a b)}{1-(a b)^t}&\frac{2(a c^x)^{r t}}{\bigl(1+(a b)^t\bigr)^r}{\rm Li}_{k_1,\dots,k_{r-1},k_r-1}\bigl(1-(a b)^{-t}\bigr)\\
&=\frac{r(a b)^t}{1+(a b)^t}\log(a b)\sum_{n=0}^\infty E_n^{(k_1,\dots,k_r)}(x)\frac{t^n}{n!}+\sum_{n=0}^\infty E_{n+1}\frac{t^n}{n!}\,. 
\end{align*} 
Since 
\begin{align*}  
\frac{\log(a b)}{1-(a b)^t}&=\log(a b)\sum_{j=0}^\infty(a b)^{j t}=\log(a b)\sum_{j=0}^\infty e^{j t\log(a b)}\\
&=\log(a b)\sum_{j=0}^\infty\sum_{\nu=0}^\infty\frac{\bigl(j\log(a b)\bigr)^\nu}{\nu!}t^\nu 
\end{align*} 
and 
\begin{align*} 
\frac{r(a b)^t}{1+(a b)^t}&=r\sum_{j=0}^\infty(a b)^{-j t}\\
&=r\sum_{j=0}^\infty\sum_{\nu=0}^\infty\frac{\bigl(-j\log(a b)\bigr)^\nu}{\nu!}t^\nu\,,
\end{align*} 
we obtain 
\begin{align*} 
r\log(a c^x)\sum_{n=0}^\infty & E_n^{(k_1,\dots,k_r)}(x)\frac{t^n}{n!}\\
+\log(a b)&\sum_{j=0}^\infty\sum_{\nu=0}^\infty\bigl(j\log(a b)\bigr)^\nu\frac{t^\nu}{\nu!}\sum_{\mu=0}^\infty E_\mu^{(k_1,\dots,k_{r-1},k_r-1)}(x)\frac{t^\mu}{\mu!}\\
&=r\sum_{j=0}^\infty\sum_{\nu=0}^\infty\bigl(-j\log(a b)\bigr)^\nu\frac{t^\nu}{\nu!}\sum_{\mu=0}^\infty E_\mu\frac{t^\mu}{\mu!}+\sum_{n=0}^\infty E_{n+1}^{(k_1,\dots,k_r)}\frac{t^n}{n!}\,. 
\end{align*}  
Hence, 
\begin{align*} 
-r\log(a c^x)\sum_{n=0}^\infty &E_n^{(k_1,\dots,k_r)}(x)\frac{t^n}{n!}+\sum_{n=0}^\infty E_{n+1}^{(k_1,\dots,k_r)}(x)\frac{t^n}{n!}\\
&=r\sum_{n=0}^\infty\left(\sum_{\nu=0}^\infty\binom{n}{\nu}\sum_{j=0}^\infty\bigl(-j\log(a b)\bigr)^{n-\nu}E_\nu^{(k_1,\dots,k_r)}(x)\right)\frac{t^n}{n!}\\
&\;\;\;\;-\log(a b)\sum_{n=0}^\infty\left(\sum_{\nu=0}^\infty\binom{n}{\nu}\sum_{j=0}^\infty\bigl(j\log(a b)\bigr)^{n-\nu}E_\nu^{(k_1,\dots,k_{r-1},k_r-1)}(x)\right)\frac{t^n}{n!}\,. 
\end{align*} 
Comparing the coefficients on both sides, we have 
\begin{align*} 
E_{n+1}^{(k_1,\dots,k_r)}(x)&-r\log(a c^x)E_n^{(k_1,\dots,k_r)}(x)\\
&=r\sum_{\nu=0}^\infty\binom{n}{\nu}\sum_{j=0}^\infty\bigl(-j\log(a b)\bigr)^{n-\nu}E_\nu^{(k_1,\dots,k_r)}(x)\\
&\;\;\;\;-\log(a b)\sum_{\nu=0}^\infty\binom{n}{\nu}\sum_{j=0}^\infty\bigl(j\log(a b)\bigr)^{n-\nu}E_\nu^{(k_1,\dots,k_{r-1},k_r-1)}(x)\,. 
\end{align*} 

If $k_r=1$, then the second term on the left-hand side becomes 
$$
2(a c^x)^{r t}\frac{-(a b)^{-t}\log(a b)}{(a b)^{-t}}{\rm Li}_{k_1,\dots,k_{r-1}}\,. 
$$ 
After dividing of $\bigl(1+(a b)^t\bigr)^r$, this second term becomes 
\begin{align*}  
-\log(a b)&\frac{2(a c^x)^{r t}}{\bigl(1+(a b)^t\bigr)^r}{\rm Li}_{k_1,\dots,k_{r-1}}\\
&=-\log(a b)\sum_{n=0}^\infty E_n^{(k_1,\dots,k_{r-1})}(x)\frac{t^n}{n!}\,. 
\end{align*} 
Finally, we get 
\begin{align*} 
E_{n+1}^{(k_1,\dots,k_r)}(x)&-r\log(a c^x)E_n^{(k_1,\dots,k_r)}(x)\\
&=r\sum_{\nu=0}^\infty\binom{n}{\nu}\sum_{j=0}^\infty\bigl(-j\log(a b)\bigr)^{n-\nu}E_\nu^{(k_1,\dots,k_r)}(x)
+\log(a b)E_n^{(k_1,\dots,k_{r-1})}(x)\,. 
\end{align*}
\end{proof}
\end{thm}

\section{Symmetrized Generalization of ${E}^{(k_1, k_2,\ldots, k_r)}_n(x;a,b,c)$}
The poly-Euler polynomials ${E}^{(-k)}_n(x;a,b,c)$ have been given symmetrized generalization \cite{CorJol1} of the form
\begin{equation*}
D^{(m)}_n(x,y;a,b,c)=\frac{1}{(\ln a+\ln b)^n}\sum_{k=0}^m\binom{m}{k}{E}^{(-k)}_n(x;a,b,c)\left(\frac{y\ln c+\ln a}{\ln a +\ln b}\right)^{m-k}.
\end{equation*}
This symmetrized generalization possesses the following double generating function 
\begin{equation}\label{symm1}
\sum_{n=0}^{\infty}\sum_{m=0}^{\infty}D^{(m)}_n(x,y;a,b,c)\frac{t^n}{n!}\frac{u^m}{m!}=\frac{2e^{\left(\frac{y\ln c+\ln a}{\ln a +\ln b}\right)u}e^{\left(\frac{x\ln c+\ln a}{\ln a +\ln b}\right)t}e^{t+u}\left(1-e^{-t}\right)}{(e^t+1)(e^t+e^u-e^{t+u})}
\end{equation}
and explicit formula
\begin{align}
D^{(m)}_n(x,y;a,b,c)&=2\sum_{j=0}^{\infty}(j!)^2\left(\sum_{l=0}^n\sum_{i=0}^{\infty}(-1)^i\frac{\left(\ln c^{x}a^{i+2}b^{i+1}\right)^{n-l}-\left(\ln c^{x}a^{i+1}b^{i}\right)^{n-l}}{(\ln a +\ln b)^{n-l}}\binom{n}{l}\braced{l}{j}\right)\times\nonumber\\
&\;\;\;\;\;\times\left(\sum_{r=0}^m\left(\frac{y\ln c+2\ln a+\ln b}{\ln a +\ln b}\right)^{m-r}\binom{m}{r}\braced{r}{j}\right).\label{symm2}
\end{align}

\begin{defn}\label{defn5}
For $m, n\ge0$, we define
\begin{equation}\label{symgen}
\mathcal{D}^{(m)}_n(x,y;a,b,c)=\sum_{k_1+k_2+\ldots +k_r=m}\binom{m}{k_1,k_2,\ldots k_r}\frac{{E}^{(-k_1,-k_2,\ldots -k_{r-1})}_n(x;a,b,c)}{(\ln a+\ln b)^n}\left(\frac{(r-1)y\ln c+\ln a}{\ln a +\ln b}\right)^{k_r}.
\end{equation}
\end{defn}

The following theorem contains the double generating function for $D^{(m)}_n(x,y;a,b,c)$.

\bigskip
\begin{thm}\label{thmm4}
For $n,m\ge0$, we have
\begin{equation}\label{eqnnnew6}
\sum_{n=0}^{\infty}\sum_{m=0}^{\infty}\mathcal{D}^{(m)}_n(x,y;a,b,c)\frac{t^n}{n!}\frac{u^m}{m!}=\frac{2e^{\left(\frac{(r-1)y\ln c+\ln a}{\ln a +\ln b}\right)u}e^{(r-1)\left(\frac{(r-1)x\ln c+\ln a}{\ln a +\ln b}\right)t}e^{\binom{r}{2}u+(r-1)t}(1-e^{-t})^{r-1}}{(1+e^{t})^{r-1}\prod_{i=1}^{r-1}(e^t+e^{iu}-e^{t+iu})}.
\end{equation}
\begin{proof}
$$\sum_{n=0}^{\infty}\sum_{m=0}^{\infty}\mathcal{D}^{(m)}_n(x,y;a,b,c)\frac{t^n}{n!}\frac{u^m}{m!}\qquad\qquad\qquad\qquad\qquad\qquad\qquad\qquad\qquad\qquad\qquad\qquad\qquad\qquad\qquad\qquad$$
\begin{eqnarray*}
&=&\sum_{n=0}^{\infty}\sum_{m=0}^{\infty}\sum_{k_1+k_2+\ldots +k_r=m}\frac{{E}^{(-k_1,-k_2,\ldots -k_{r-1})}_n(x;a,b,c)}{(\ln a+\ln b)^n}\left(\frac{(r-1)y\ln c+\ln a}{\ln a +\ln b}\right)^{k_r}\frac{t^n}{n!}\times\\
&&\;\;\;\;\;\;\times\frac{u^m}{k_1!k_2!\ldots k_r!}\\
&=&\sum_{n=0}^{\infty}\sum_{k_1+k_2+\ldots +k_r\ge0}\frac{{E}^{(-k_1,-k_2,\ldots -k_{r-1})}_n(x;a,b,c)}{(\ln a+\ln b)^n}\left(\frac{(r-1)y\ln c+\ln a}{\ln a +\ln b}\right)^{k_r}\frac{t^n}{n!}\times\\
&&\;\;\;\;\;\;\times\frac{u^{k_1+k_2+\ldots +k_r}}{k_1!k_2!\ldots k_r!}\\
&=&\sum_{n=0}^{\infty}\sum_{k_1+k_2+\ldots +k_{r-1}\ge0}\frac{{E}^{(-k_1,-k_2,\ldots -k_{r-1})}_n(x;a,b,c)}{(\ln a+\ln b)^n}\sum_{k_r\ge0}\left(\frac{(r-1)y\ln c+\ln a}{\ln a +\ln b}\right)^{k_r}\frac{u^{k_r}}{k_r!}\times\\
&&\;\;\;\;\;\;\times\frac{t^n}{n!}\frac{u^{k_1+k_2+\ldots +k_{r-1}}}{k_1!k_2!\ldots k_{r-1}!}\\
&=&e^{\left(\frac{(r-1)y\ln c+\ln a}{\ln a +\ln b}\right)u}\sum_{n=0}^{\infty}\sum_{k_1+k_2+\ldots +k_{r-1}\ge0}\frac{{E}^{(-k_1,-k_2,\ldots -k_{r-1})}_n(x;a,b,c)}{(\ln a+\ln b)^n}\frac{t^n}{n!}\frac{u^{k_1+k_2+\ldots +k_{r-1}}}{k_1!k_2!\ldots k_{r-1}!}
\end{eqnarray*}
Using identity (\ref{prop2}), we obtain
$$\sum_{n=0}^{\infty}\sum_{m=0}^{\infty}\mathcal{D}^{(m)}_n(x,y;a,b,c)\frac{t^n}{n!}\frac{u^m}{m!}\qquad\qquad\qquad\qquad\qquad\qquad\qquad\qquad\qquad\qquad\qquad\qquad\qquad\qquad\qquad\qquad$$
\begin{eqnarray*}
&=&e^{\left(\frac{(r-1)y\ln c+\ln a}{\ln a +\ln b}\right)u}\sum_{k_1+k_2+\ldots +k_{r-1}\ge0}\sum_{n=0}^{\infty}{E}^{(-k_1,-k_2,\ldots -k_{r-1})}_n\left(\frac{(r-1)x\ln c+\ln a}{\ln a +\ln b}\right)\frac{t^n}{n!}\times\\
&&\;\;\;\;\;\;\times\frac{u^{k_1+k_2+\ldots +k_{r-1}}}{k_1!k_2!\ldots k_{r-1}!}\\
&=&e^{\left(\frac{(r-1)y\ln c+\ln a}{\ln a +\ln b}\right)u}e^{(r-1)\left(\frac{(r-1)x\ln c+\ln a}{\ln a +\ln b}\right)t}\sum_{k_1+k_2+\ldots +k_{r-1}\ge0}\frac{2{\rm Li}_{(-k_1, -k_2,\ldots, -k_{r-1})}(1-e^{-t})}{(1+e^{t})^{r-1}}\times\\
&&\;\;\;\;\;\;\times\frac{u^{k_1+k_2+\ldots +k_{r-1}}}{k_1!k_2!\ldots k_{r-1}!}\\
&=&\frac{2e^{\left(\frac{(r-1)y\ln c+\ln a}{\ln a +\ln b}\right)u}e^{(r-1)\left(\frac{(r-1)x\ln c+\ln a}{\ln a +\ln b}\right)t}}{(1+e^{t})^{r-1}}\sum_{0<m_1<m_2<\ldots <m_{r-1}}(1-e^{-t})^{m_{r-1}}\mathcal{S}(u,m_1,m_2,\ldots, m_{r-1})
\end{eqnarray*}
where
$$\mathcal{S}(u,m_1,m_2,\ldots, m_{r-1})=\sum_{k_1+k_2+\ldots +k_{r-1}\ge0}\frac{(um_1)^{k_1}(um_2)^{k_2}\ldots (um_{r-1})^{k_{r-1}}}{k_1!k_2!\ldots k_{r-1}!}\qquad\qquad\qquad\qquad\qquad\qquad\qquad\qquad\qquad$$
\begin{eqnarray*}
&=&\sum_{\widehat{m}\ge0}\frac{1}{\widehat{m}!}\sum_{k_1+k_2+\ldots +k_{r-1}=\widehat{m}}\binom{\widehat{m}}{k_1, k_2, \ldots k_{r-1}}(um_1)^{k_1}(um_2)^{k_2}\ldots (um_{r-1})^{k_{r-1}}\\
&=&\sum_{\widehat{m}\ge0}\frac{(um_1+um_2+\ldots +um_{r-1})^{\widehat{m}}}{\widehat{m}!}\\
&=&e^{u(m_1+m_2+\ldots +m_{r-1})}.
\end{eqnarray*}
Thus,
$$\sum_{n=0}^{\infty}\sum_{m=0}^{\infty}\mathcal{D}^{(m)}_n(x,y;a,b,c)\frac{t^n}{n!}\frac{u^m}{m!}\qquad\qquad\qquad\qquad\qquad\qquad\qquad\qquad\qquad\qquad\qquad\qquad\qquad\qquad\qquad\qquad$$
\begin{eqnarray*}
&=&\frac{2e^{\left(\frac{(r-1)y\ln c+\ln a}{\ln a +\ln b}\right)u}e^{(r-1)\left(\frac{(r-1)x\ln c+\ln a}{\ln a +\ln b}\right)t}}{(1+e^{t})^{r-1}}\sum_{0<m_1<m_2<\ldots <m_{r-1}}(1-e^{-t})^{m_{r-1}}e^{u(m_1+m_2+\ldots +m_{r-1})}\\
&=&\frac{2e^{\left(\frac{(r-1)y\ln c+\ln a}{\ln a +\ln b}\right)u}e^{(r-1)\left(\frac{(r-1)x\ln c+\ln a}{\ln a +\ln b}\right)t}}{(1+e^{t})^{r-1}}\times\\
&&\;\;\;\;\;\;\times\frac{e^u(1-e^{-t})}{1-e^u(1-e^{-t})}\frac{e^{2u}(1-e^{-t})}{1-e^{2u}(1-e^{-t})}\ldots\frac{e^{(r-1)u}(1-e^{-t})}{1-e^{(r-1)u}(1-e^{-t})}\\
&=&\frac{2e^{\left(\frac{(r-1)y\ln c+\ln a}{\ln a +\ln b}\right)u}e^{(r-1)\left(\frac{(r-1)x\ln c+\ln a}{\ln a +\ln b}\right)t}e^{\binom{r}{2}u}(1-e^{-t})^{r-1}}{(1+e^{t})^{r-1}\prod_{i=1}^{r-1}(1-e^{iu}(1-e^{-t}))}\\
&=&\frac{2e^{\left(\frac{(r-1)y\ln c+\ln a}{\ln a +\ln b}\right)u}e^{(r-1)\left(\frac{(r-1)x\ln c+\ln a}{\ln a +\ln b}\right)t}e^{\binom{r}{2}u+(r-1)t}(1-e^{-t})^{r-1}}{(1+e^{t})^{r-1}\prod_{i=1}^{r-1}(e^t+e^{iu}-e^{t+iu})}.
\end{eqnarray*}
\end{proof}
\end{thm}

\smallskip
Note that equation (\ref{symm1}) can easily be deduced from equation (\ref{eqnnnew6}) by taking $r=1$. It is then interesting to establish an explicit formula for $\mathcal{D}^{(m)}_n(x,y;a,b,c)$ parallel to equation (\ref{symm2}). To do this, let us consider first the following expression from the right-hand side of equation (\ref{eqnnnew6}). That is,
\begin{eqnarray*}
\frac{1}{(1+e^{t})^{r-1}\prod_{i=1}^{r-1}(e^t+e^{iu}-e^{t+iu})}&=&\left(\sum_{n=0}^{\infty}(-e)^{nt}\right)^{r-1}\prod_{i=1}^{r-1}\frac{1}{(1- (e^t-1)(e^{iu}-1))}\\
&=&\left(\sum_{n=0}^{\infty}(-e)^{nt}\right)^{r-1}\prod_{i=1}^{r-1}\sum_{j=0}^{\infty}(e^t-1)^j(e^{iu}-1)^j\\
&=&\left(\sum_{n=0}^{\infty}(-e)^{nt}\right)^{r-1}\prod_{i=1}^{r-1}\sum_{c_i\ge0}(e^t-1)^{c_i}(e^{iu}-1)^{c_i}\\
\left(\sum_{n=0}^{\infty}(-e)^{nt}\right)^{r-1}&=&\sum_{q=0}^{\infty}\sum_{k_{r-1}=0}^n\sum_{k_{r-2}=0}^{n-k_{r-1}}\ldots\sum_{k_{1}=0}^{n-k_{r-1}-\ldots -k_{2}}(-1)^qe^{qt}\\
&=&\sum_{q=0}^{\infty}(-1)^q\frac{\prod_{j=0}^{q-2}(q+1+j)}{(q-1)!}e^{qt}
\end{eqnarray*}
\begin{eqnarray*}
\prod_{i=1}^{r-1}\sum_{c_i\ge0}(e^t-1)^{c_i}(e^{iu}-1)^{c_i}&=&\sum_{j=0}^{\infty}\sum_{c_1+c_2+\ldots +c_{r-1}=j}(e^t-1)^{j}\prod_{i=1}^{r-1}(e^{iu}-1)^{c_i}\\
&=&\sum_{j=0}^{\infty}\sum_{c_1+c_2+\ldots +c_{r-1}=j}j!\sum_{n=0}^{\infty}\braced{n}{j}\frac{t^n}{n!}\prod_{i=1}^{r-1}c_i!\sum_{n=0}^{\infty}\braced{n}{c_i}\frac{u^n}{n!}
\end{eqnarray*}
\begin{eqnarray*}
\prod_{i=1}^{r-1}c_i!\sum_{n=0}^{\infty}\braced{n}{c_i}\frac{u^n}{n!}&=&\sum_{m=0}^{\infty}\sum_{d_{r-2}=0}^{m}\sum_{d_{r-3}=0}^{m-d_{r-2}}\ldots \sum_{d_{1}=0}^{m-d_{r-2}-\ldots - d_2}c_1!\braced{m-d_{r-2}-\ldots - d_{1}}{c_1}\times\\
&&\;\;\;\times\prod_{i=1}^{r-2}\binom{m-d_{r-2}-\ldots - d_{i+1}}{d_i}\braced{d_i}{c_{i+1}}c_{i+1}!(i+1)^{d_i}\frac{u^m}{m!}
\end{eqnarray*}
$$e^{(r-1)\left(\frac{(r-1)x\ln c+\ln a}{\ln a +\ln b}\right)t}(1-e^{-t})^{r-1}\left(\sum_{n=0}^{\infty}(-e)^{nt}\right)^{r-1}\left(j!\sum_{n=0}^{\infty}\braced{n}{j}\frac{t^n}{n!}\right)
\qquad\qquad\qquad\qquad\qquad\qquad\qquad$$
\begin{eqnarray*}
&=&e^{(r-1)\left(\frac{(r-1)x\ln c+\ln a}{\ln a +\ln b}\right)t}\left(\sum_{k=0}^{r-1}(-1)^k\binom{r-1}{k}e^{-kt}\right)\left(\sum_{q=0}^{\infty}(-1)^q\frac{\prod_{j=0}^{q-2}(q+1+j)}{(q-1)!}e^{qt}\right)\times\\
&&\;\;\;\;\;\;\times\left(j!\sum_{n=0}^{\infty}\braced{n}{j}\frac{t^n}{n!}\right)\\
&=&\left(\sum_{q=0}^{\infty}\sum_{k=0}^{r-1}e^{\left(\frac{(r-1)^2x\ln c+(q-k)\ln b+(q-k+r-1)\ln a}{\ln a+\ln b}\right)}\frac{(-1)^{k+q}\binom{r-1}{k}\prod_{j=0}^{q-2}(q+1+j)}{(q-1)!}\right)\times\\
&&\;\;\;\;\;\;\times\left(j!\sum_{n=0}^{\infty}\braced{n}{j}\frac{t^n}{n!}\right)
\end{eqnarray*}
\begin{eqnarray*}
&=&\left(\sum_{n=0}^{\infty}\sum_{q=0}^{\infty}\sum_{k=0}^{r-1}\left(\frac{(r-1)^2x\ln c+(q-k)\ln b+(q-k+r-1)\ln a}{\ln a+\ln b}\right)^{n}\right.\\
&&\;\;\;\;\;\left.\frac{(-1)^{k+q}\binom{r-1}{k}\prod_{j=0}^{q-2}(q+1+j)}{(q-1)!}\frac{t^n}{n!}\right)\left(j!\sum_{n=0}^{\infty}\braced{n}{j}\frac{t^n}{n!}\right)\\
&=&\sum_{n=0}^{\infty}\left(\sum_{p=0}^n\sum_{q=0}^{\infty}\sum_{k=0}^{r-1}\binom{n}{p}\left(\frac{(r-1)^2x\ln c+(q-k)\ln b+(q-k+r-1)\ln a}{\ln a+\ln b}\right)^{n-p}\right.\\
&&\;\;\;\;\;\left.\frac{(-1)^{k+q}\binom{r-1}{k}\prod_{j=0}^{q-2}(q+1+j)}{(q-1)!}j!\braced{p}{j}\right)\frac{t^n}{n!}
\end{eqnarray*}
$$e^{\left(\frac{(r-1)y\ln c+\ln a}{\ln a +\ln b}\right)u}e^{\binom{r}{2}u}\left(\prod_{i=1}^{r-1}c_i!\sum_{n=0}^{\infty}\braced{n}{c_i}\frac{u^n}{n!}\right)
\qquad\qquad\qquad\qquad\qquad\qquad\qquad\qquad\qquad\qquad\qquad$$
\begin{eqnarray*}
&=&\left(\sum_{m=0}^{\infty}\left(\frac{(r-1)y\ln c+\binom{r}{2}\ln b+\left\{\binom{r}{2}+1\right\}\ln a}{\ln a +\ln b}\right)^m\frac{u^m}{m!}\right)\left(\sum_{m=0}^{\infty}\sum_{d_{r-2}=0}^{m}\sum_{d_{r-3}=0}^{m-d_{r-2}}\ldots\right.\\
&&\;\;\;\left. \sum_{d_{1}=0}^{m-d_{r-2}-\ldots - d_2}c_1!\braced{m-d_{r-2}-\ldots - d_{1}}{c_1}\prod_{i=1}^{r-2}\binom{m-d_{r-2}-\ldots - d_{i+1}}{d_i}\times\right.\\
&&\;\;\;\;\;\;\times\left. \braced{d_i}{c_{i+1}}c_{i+1}!(i+1)^{d_i}\frac{u^m}{m!}\right)\\
&=&\sum_{m=0}^{\infty}\left\{\sum_{l=0}^m\left(\frac{(r-1)y\ln c+\binom{r}{2}\ln b+\left\{\binom{r}{2}+1\right\}\ln a}{\ln a +\ln b}\right)^{m-l}\sum_{d_{r-2}=0}^{l}\sum_{d_{r-3}=0}^{l-d_{r-2}}\ldots\right.\\
&&\;\;\;\left. \sum_{d_{1}=0}^{l-d_{r-2}-\ldots - d_2}c_1!\braced{l-d_{r-2}-\ldots - d_{1}}{c_1}\prod_{i=1}^{r-2}\binom{l-d_{r-2}-\ldots - d_{i+1}}{d_i}\times\right.\\
&&\;\;\;\;\;\;\times\left. \braced{d_i}{c_{i+1}}c_{i+1}!(i+1)^{d_i}\right\}\frac{u^m}{m!}
\end{eqnarray*}
$$\sum_{n=0}^{\infty}\sum_{m=0}^{\infty}\mathcal{D}^{(m)}_n(x,y;a,b,c)\frac{t^n}{n!}\frac{u^m}{m!}\qquad\qquad\qquad\qquad\qquad\qquad\qquad\qquad\qquad\qquad\qquad\qquad\qquad\qquad\qquad\qquad$$
\begin{eqnarray*}
&=&\sum_{n=0}^{\infty}\sum_{m=0}^{\infty}\left\{\sum_{j=0}^{\infty}\sum_{c_1+c_2+\ldots +c_{r-1}=j}\sum_{l=0}^m\left(\frac{(r-1)y\ln c+\binom{r}{2}\ln b+\left\{\binom{r}{2}+1\right\}\ln a}{\ln a +\ln b}\right)^{m-l}\right.\\
&&\left. \sum_{d_{r-2}=0}^{l}\sum_{d_{r-3}=0}^{l-d_{r-2}}\ldots\sum_{d_{1}=0}^{l-d_{r-2}-\ldots - d_2}c_1!\braced{l-d_{r-2}-\ldots - d_{1}}{c_1}\prod_{i=1}^{r-2}\binom{l-d_{r-2}-\ldots - d_{i+1}}{d_i}\times\right.\\
&&\left. \times\braced{d_i}{c_{i+1}}c_{i+1}!(i+1)^{d_i}\sum_{p=0}^n\sum_{q=0}^{\infty}\sum_{k=0}^{r-1}\binom{n}{p}\frac{(-1)^{k+q}\binom{r-1}{k}\prod_{j=0}^{q-2}(q+1+j)}{(q-1)!}j!\braced{p}{j}F_q^p\right\}\frac{t^n}{n!}\frac{u^m}{m!}
\end{eqnarray*}
where 
$$F_q^p=\left(\frac{(r-1)^2x\ln c+(q-k)\ln b+(q-k+r-1)\ln a}{\ln a+\ln b}\right)^{n-p}.$$
Comparing coefficients, we obtain the following theorem.

\bigskip
\begin{thm}\label{thm555}
For $n,m\ge0$, we have
$$\mathcal{D}^{(m)}_n(x,y;a,b,c)=\sum_{j=0}^{\infty}\sum_{c_1+c_2+\ldots +c_{r-1}=j}\sum_{l=0}^m\left(\frac{(r-1)y\ln c+\binom{r}{2}\ln b+\left\{\binom{r}{2}+1\right\}\ln a}{\ln a +\ln b}\right)^{m-l}\times$$
\begin{eqnarray*}
&&\;\times\sum_{d_{r-2}=0}^{l}\sum_{d_{r-3}=0}^{l-d_{r-2}}\ldots\sum_{d_{1}=0}^{l-d_{r-2}-\ldots - d_2}c_1!\braced{l-d_{r-2}-\ldots - d_{1}}{c_1}\prod_{i=1}^{r-2}\binom{l-d_{r-2}-\ldots - d_{i+1}}{d_i}\times\\
&&\;\times\braced{d_i}{c_{i+1}}c_{i+1}!(i+1)^{d_i}\sum_{p=0}^n\sum_{q=0}^{\infty}\sum_{k=0}^{r-1}\binom{n}{p}\frac{(-1)^{k+q}\binom{r-1}{k}\prod_{j=0}^{q-2}(q+1+j)}{(q-1)!}j!\braced{p}{j}F_q^p
\end{eqnarray*}
where 
$$F_q^p=\left(\frac{(r-1)^2x\ln c+(q-k)\ln b+(q-k+r-1)\ln a}{\ln a+\ln b}\right)^{n-p}.$$
\end{thm}

\section{Generalized Multi Poly-Bernoulli Polynomials}
Parallel to the definition of generalized multi poly-Euler polynomials in (\ref{multipolyeuler}), we have the following generalization of poly-Bernoulli numbers.

\begin{defn}\label{multipolybernoulli}\rm
The generalized multi poly-Bernoulli polynomials are defined by
\begin{equation}\label{multipolybern}
\frac{{\rm Li}_{(k_1, k_2,\ldots, k_r)}(1-(ab)^{-t})}{(b^{t}-a^{-t})^r}c^{rxt}=\sum_{n=0}^{\infty}{B}^{(k_1, k_2,\ldots, k_r)}_n(x;a,b,c)\frac{t^n}{n!}.
\end{equation}
\end{defn}

One can easily prove the following theorem using the same argument in deriving the identities in Theorem \ref{add2}-\ref{add5}. 

\begin{thm}\label{addd} The generalized multi poly-Bernoulli polynomials satisfy the following identities.
\begin{align*}
{B}^{(k_1, k_2,\ldots, k_r)}_n(x;a,b,c)&=\sum_{m=0}^{\infty}\sum_{l=m}^n(r\log c)^l\braced{l}{m}\binom{n}{l}{B}^{(k_1, k_2,\ldots, k_r)}_{n-l}(-m\log c;a,b)(x)^{(m)}\\
{B}^{(k_1, k_2,\ldots, k_r)}_n(x;a,b,c)&=\sum_{m=0}^{\infty}\sum_{l=m}^n(r\log c)^l\braced{l}{m}\binom{n}{l}{B}^{(k_1, k_2,\ldots, k_r)}_{n-l}(0;a,b)(x)^{(m)}\\
{B}^{(k_1, k_2,\ldots, k_r)}_n(x;a,b,c)&=\sum_{m=0}^{\infty}\binom{n}{m}\sum_{l=0}^{n-m}\frac{\binom{n-m}{l}}{\binom{l+s}{l}}\braced{l+s}{s}{B}^{(k_1, k_2,\ldots, k_r)}_{n-m-l}(0;a,b)B^{(s)}_m(xr\log c)\\
{B}^{(k_1, k_2,\ldots, k_r)}_n(x;a,b,c)&=\sum_{m=0}^{n}\frac{\binom{n}{m}}{(1-\lambda)^s}\sum_{j=0}^s\binom{s}{j}(-\lambda)^{s-j}{B}^{(k_1, k_2,\ldots, k_r)}_{n-m}(j;a,b)H^{(s)}_m(xr\log c;\lambda)
\end{align*}

\end{thm}

\smallskip
The next theorem contains an explicit formula for ${B}^{(k_1, k_2,\ldots, k_r)}_n (x; a, b, c)$.

\bigskip
\begin{thm}\label{thm40} {\rm ({\bf Explicit Formula})} For $k \in \mathbb{Z}$, $n\geq0$, we have
\begin{equation}
{B}^{(k_1, k_2,\ldots, k_r)}_n (x; a, b, c) =\sum_{m_r>\ldots >m_1>0}\frac{1}{m_1^{k_1}m_2^{k_2}\ldots m_r^{k_r}}\sum_{j=0}^{m_r-r}(-1)^j\binom{m_r-r}{j}(rx-j\ln a-(j+1)\ln b)^n.
\end{equation}
\begin{proof}
\begin{align*}
\frac{{\rm Li}_{(k_1, k_2,\ldots, k_r)}(1-(ab)^{-t})}{(b^t-a^{-t})^r}&=b^{-rt}\left(\sum_{m_r>\ldots >m_1>0}\frac{(1-(ab)^{-t})^{m_r-r}}{m_1^{k_1}m_2^{k_2}\ldots m_r^{k_r}}\right)\\
&=b^{-rt}\sum_{m_r>\ldots >m_1>0}\frac{1}{m_1^{k_1}m_2^{k_2}\ldots m_r^{k_r}}\sum_{j=0}^{m_r-r}(-1)^j\binom{m_r-r}{j}e^{-jt\ln (ab)}\\
&=\sum_{m_r>\ldots >m_1>0}\frac{1}{m_1^{k_1}m_2^{k_2}\ldots m_r^{k_r}}\sum_{j=0}^{m_r-r}(-1)^j\binom{m_r-r}{j}e^{-t(j\ln a+(j+1)\ln b)}.
\end{align*}
So, we get

\smallskip
\begin{equation*}
\frac{{\rm Li}_{(k_1, k_2,\ldots, k_r)}(1-(ab)^{-t})}{(b^t-a^{-t})^r}e^{xrt\ln c}\qquad\qquad\qquad\qquad\qquad\qquad\qquad\qquad\qquad\qquad\qquad\qquad\qquad\qquad\qquad\qquad\qquad\qquad
\end{equation*}
\begin{align*}
&=\sum_{m_r>\ldots >m_1>0}\frac{1}{m_1^{k_1}m_2^{k_2}\ldots m_r^{k_r}}\sum_{j=0}^{m_r-r}(-1)^j\binom{m_r-r}{j}e^{t(rx-j\ln a-(j+1)\ln b)}\\
&=\sum_{n=0}^{\infty}\left(\sum_{m_r>\ldots >m_1>0}\frac{1}{m_1^{k_1}m_2^{k_2}\ldots m_r^{k_r}}\sum_{j=0}^{m_r-r}(-1)^j\binom{m_r-r}{j}(rx-j\ln a-(j+1)\ln b)^n\right)\frac{t^n}{n!}
\end{align*}
By comparing the coefficients of $\frac{t^n}{n!}$ on both sides, the proof is completed.
\end{proof}
\end{thm}

\bigskip
The next theorem contains an expression of ${B}^{(k_1, k_2,\ldots, k_r)}_n(x;a,b,c)$ as polynomial in $x$.

\bigskip
\begin{thm}\label{thm41}
The generalized multi poly-Bernoulli polynomials satisfy the following relation
\begin{equation}\label{eqnnew1}
{B}^{(k_1, k_2,\ldots, k_r)}_n(x;a,b,c)=\sum_{i=0}^n\binom{n}{i}(\ln c)^{n-i}{B}^{(k_1, k_2,\ldots, k_r)}_i(a,b)x^{n-i}
\end{equation}
\begin{proof}
Using (\ref{multipolybern}), we have
\begin{eqnarray*}
\sum_{n=0}^{\infty}{B}^{(k_1, k_2,\ldots, k_r)}_n(x;a,b,c)\frac{t^n}{n!}&=&\frac{{\rm Li}_{(k_1, k_2,\ldots, k_r)}(1-(ab)^{-t})}{(b^t-a^{-t})^r}c^{xt}=e^{xt\ln c}\sum_{n=0}^{\infty}{B}^{(k_1, k_2,\ldots, k_r)}_n(a,b)\frac{t^n}{n!}\\
&=&\sum_{n=0}^{\infty}\sum_{i=0}^n\frac{(xt\ln c)^{n-i}}{(n-i)!}{B}^{(k_1, k_2,\ldots, k_r)}_i(a,b)\frac{t^{i}}{i!}\\
&=&\sum_{n=0}^{\infty}\left(\sum_{i=0}^n\binom{n}{i}(\ln c)^{n-i}{B}^{(k_1, k_2,\ldots, k_r)}_i(a,b)x^{n-i}\right)\frac{t^{n}}{n!}.
\end{eqnarray*}
Comparing the coefficients of $\frac{t^{n}}{n!}$, we obtain the desired result.
\end{proof}
\end{thm}

Note that, when $a=c=e$ and $b=1$, Definition \ref{multipolybernoulli} reduces to
\begin{equation}\label{multipolybernpolynomial}
\frac{{\rm Li}_{(k_1, k_2,\ldots, k_r)}(1-e^{-t})}{(1-e^{-t})^r}e^{rxt}=\sum_{n=0}^{\infty}{B}^{(k_1, k_2,\ldots, k_r)}_n(x)\frac{t^n}{n!}.
\end{equation}
The following theorem gives a relation between ${B}^{(k_1, k_2,\ldots, k_r)}_n(x;a,b,c)$ and ${B}^{(k_1, k_2,\ldots, k_r)}_n(x)$.

\bigskip
\begin{thm}\label{thm42}
The generalized multi poly-Bernoulli polynomials satisfy the following relation
\begin{equation}\label{eqnnew2}
{B}^{(k_1, k_2,\ldots, k_r)}_n(x;a,b,c)=(\ln a+\ln b)^n{B}^{(k_1, k_2,\ldots, k_r)}_n\left(\frac{x\ln c-r\ln b}{\ln a+\ln b}\right)
\end{equation}
\begin{proof}
Using (\ref{multipolybern}), we have
\begin{eqnarray*}
\sum_{n=0}^{\infty}{B}^{(k_1, k_2,\ldots, k_r)}_n(x;a,b,c)\frac{t^n}{n!}&=&\frac{{\rm Li}_{(k_1, k_2,\ldots, k_r)}(1-(ab)^{-t})}{b^{rt}(1-(ab)^{-t})^r}e^{xt\ln c}\\
&=&e^{\frac{x\ln c-r\ln b}{\ln ab}t\ln ab}\frac{{\rm Li}_{(k_1, k_2,\ldots, k_r)}(1-e^{-t\ln ab})}{1+e^{-t\ln ab}}\\
&=&\sum_{n=0}^{\infty}(\ln a+\ln b)^n{B}^{(k_1, k_2,\ldots, k_r)}_n\left(\frac{x\ln c-r\ln b}{\ln a+\ln b}\right)\frac{t^{n}}{n!}.
\end{eqnarray*}
Comparing the coefficients of $\frac{t^{n}}{n!}$, we obtain the desired result.
\end{proof}
\end{thm}

\bigskip
\begin{thm}\label{thm43}
The generalized poly-Bernoulli polynomials satisfy the following relation
\begin{equation}\label{eqnnew3}
\frac{d}{dx}{B}^{(k_1, k_2,\ldots, k_r)}_{n+1}(x;a,b,c)=(n+1)(\ln c){B}^{(k_1, k_2,\ldots, k_r)}_{n}(x;a,b,c)
\end{equation}
\begin{proof}
Using (\ref{multipolybern}), we have
\begin{eqnarray*}
\sum_{n=0}^{\infty}\frac{d}{dx}{B}^{(k_1, k_2,\ldots, k_r)}_{n}(x;a,b,c)\frac{t^{n}}{n!}&=&\frac{t(\ln c){\rm Li}_{(k_1, k_2,\ldots, k_r)}(1-(ab)^{-t})}{(b^t-a^{-t})^r}e^{xrt\ln c}\\
\sum_{n=0}^{\infty}\frac{d}{dx}{B}^{(k_1, k_2,\ldots, k_r)}_{n}(x;a,b,c)\frac{t^{n-1}}{n!}&=&\sum_{n=0}^{\infty}(r\ln c){B}^{(k_1, k_2,\ldots, k_r)}_{n}(x;a,b,c)\frac{t^{n}}{n!}.
\end{eqnarray*}
Hence, 
\begin{eqnarray*}
\sum_{n=0}^{\infty}\frac{1}{n+1}\frac{d}{dx}{B}^{(k_1, k_2,\ldots, k_r)}_{n+1}(x;a,b,c)\frac{t^{n}}{n!}=\sum_{n=0}^{\infty}(\ln c){B}^{(k_1, k_2,\ldots, k_r)}_{n}(x;a,b,c)\frac{t^{n}}{n!}.
\end{eqnarray*}
Comparing the coefficients of $\frac{t^{n}}{n!}$, we obtain the desired result.
\end{proof}
\end{thm}

The following corollary immediately follows from Theorem \ref{thm43} by taking $c=e$. For brevity, let us denote ${B}^{(k_1, k_2,\ldots, k_r)}_{n}(x;a,b,e)$ by ${B}^{(k_1, k_2,\ldots, k_r)}_{n}(x;a,b)$.

\bigskip
\begin{cor}\label{cor1}
The generalized poly-Bernoulli polynomials are Appell polynomials in the sense that
\begin{equation}\label{eqnnew4}
\frac{d}{dx}{B}^{(k_1, k_2,\ldots, k_r)}_{n+1}(x;a,b)=(n+1){B}^{(k_1, k_2,\ldots, k_r)}_{n}(x;a,b)
\end{equation}
\end{cor}

\smallskip
Consequently, using the characterization of Appell polynomials \cite{Lee, Shohat, Toscano}, the following addition formula can easily be obtained.  

\bigskip
\begin{cor}\label{cor2}
The generalized poly-Bernoulli polynomials satisfy the following addition formula
\begin{equation}\label{eqnnew5}
{B}^{(k_1, k_2,\ldots, k_r)}_{n}(x+y;a,b)=\sum_{i=0}^n\binom{n}{i}{B}^{(k_1, k_2,\ldots, k_r)}_i(x;a,b)y^{n-i}
\end{equation}
\end{cor}

However, we can derive the addition formula for ${B}^{(k_1, k_2,\ldots, k_r)}_n(x;a,b,c)$ as follows
\begin{eqnarray*}
\sum_{n=0}^{\infty}{B}^{(k_1, k_2,\ldots, k_r)}_n(x+y;a,b,c)\frac{t^n}{n!}&=&\frac{{\rm Li}_{(k_1, k_2,\ldots, k_r)}(1-(ab)^{-t})}{(b^t-a^{-t})^r}c^{(x+y)rt}\\
&=&\frac{{\rm Li}_{(k_1, k_2,\ldots, k_r)}(1-(ab)^{-t})}{(b^t-a^{-t})^r}c^{xrt}c^{yrt}\\
&=&\left(\sum_{n=0}^{\infty}{B}^{(k_1, k_2,\ldots, k_r)}_n(x;a,b,c)\frac{t^n}{n!}\right)\left(\sum_{n=0}^{\infty}(yr\ln c)^n\frac{t^n}{n!}\right)\\
&=&\sum_{n=0}^{\infty}\left(\sum_{i=0}^n\binom{n}{i}(yr\ln c)^{n-i}{B}^{(k_1, k_2,\ldots, k_r)}_i(x;a,b,c)\right)\frac{t^n}{n!}.
\end{eqnarray*}
Comparing the coefficients of $\frac{t^n}{n!}$ yields the following result.

\bigskip
\begin{thm}\label{multithm2}
The generalized poly-Bernoulli polynomials satisfy the following addition formula
\begin{equation*}
{B}^{(k_1, k_2,\ldots, k_r)}_{n}(x+y;a,b,c)=\sum_{i=0}^n\binom{n}{i}(r\ln c)^{n-i}{B}^{(k_1, k_2,\ldots, k_r)}_{i}(x;a,b,c)y^{n-i}.
\end{equation*}
\end{thm}

\section{Hurwitz-Lerch Type Multi Poly-Bernoulli Polynomials}
Consider the case in which $x=1$, $a=e$ and $b=c=1$ for the parameters in Definition \ref{multipolybernoulli}. Then we have
\begin{equation}\label{multipolybernpolynomial-1}
\frac{{\rm Li}_{(k_1, k_2,\ldots, k_r)}(1-e^{-t})}{(1-e^{-t})^r}=\sum_{n=0}^{\infty}{B}^{(k_1, k_2,\ldots, k_r)}_n\frac{t^n}{n!}.
\end{equation}
This can be generalized using the following generalization of Hurwitz-Lerch multiple zeta values 
\begin{equation}\label{HLerch1}
{ \Phi}_{(k_1,k_2,\ldots, k_r)}(z,a)=\sum_{ 0\le m_1\le m_2\le\ldots \le m_r }\frac{z^{m_r}}{(m_1+a-r+1)^{k_1} (m_2+a-r+2)^{k_2}\ldots (m_r+a)^{k_r}}.
\end{equation}
Note that
\begin{align*}
{\rm Li}_{(k_1,k_2,\ldots, k_r)}(z)&=\sum_{ 0< m_1< m_2<\ldots < m_r }\frac{z^{m_r}}{m_1^{k_1} m_2^{k_2}\ldots m_r^{k_r}}\\
&=z^r\sum_{ 0< m_1< m_2<\ldots < m_r }\frac{z^{m_r-r}}{{m_1}^{k_1} {m_2}^{k_2}\ldots {m_r}^{k_r}}\\
&=z^r\sum_{ 0\le m_1\le m_2\le\ldots \le m_r }\frac{z^{m_r}}{{(m_1+1)}^{k_1} {(m_2+2)}^{k_2}\ldots {(m_r+r)}^{k_r}}\\
&=z^r{ \Phi}_{(k_1,k_2,\ldots, k_r)}(z,r)
\end{align*}
Thus, we have
\begin{equation}
\frac{{\rm Li}_{(k_1,k_2,\ldots, k_r)}(z)}{z^r}={ \Phi}_{(k_1,k_2,\ldots, k_r)}(z,r).
\end{equation}
More precisely, one can generalize (\ref{HpolyB}) as follows
\begin{equation}\label{HpolyB}
{ \Phi}_{(k_1,k_2,\ldots, k_r)}(1-e^{-t},a)=\sum_{n=0}^{\infty}B_{n,a}^{(k_1,k_2,\ldots, k_r)}\frac{t^n}{n!}.
\end{equation}
We call $B_{n,a}^{(k_1,k_2,\ldots, k_r)}$ as {\it Hurwitz-Lerch Type Multi Poly-Bernoulli Numbers}. Furthermore, we can define the {\it Hurwitz-Lerch Type Multi Poly-Bernoulli Polynomials}, denoted by $B_{n,a}^{(k_1,k_2,\ldots, k_r)}(x)$, as follows
\begin{equation}\label{HpolyB1}
{ \Phi}_{(k_1,k_2,\ldots, k_r)}(1-e^{-t},a)e^{rxt}=\sum_{n=0}^{\infty}B_{n,a}^{(k_1,k_2,\ldots, k_r)}(x)\frac{t^n}{n!}
\end{equation}
where $B_{n,a}^{(k_1,k_2,\ldots, k_r)}(0)=B_{n,a}^{(k_1,k_2,\ldots, k_r)}$.

\smallskip
The next theorem contains an explicit formula for $B_{n,a}^{(k_1,k_2,\ldots, k_r)}(x)$ expressed in terms of the $(r,\beta)$-Stirling numbers $\braced{n}{m}_{\beta,r}$ \cite{CORCA}, which satisfy the following exponential generating function
\begin{equation}\label{rbStirling}
\frac{1}{\beta^mm!}e^{rt}(e^{\beta t}-1)^m=\sum_{n=0}^{\infty}\braced{n}{m}_{\beta,r}\frac{t^n}{n!}.
\end{equation}

\bigskip
\begin{thm}\label{res1}
The Hurwitz-Lerch type multi poly-Bernoulli polynomials have the following explicit formula
\begin{equation}
B_{n,a}^{(k_1,k_2,\ldots, k_r)}(x)=\sum_{0\le m_1\le m_2\le\ldots \le m_r\le n}\frac{m_r!\braced{n}{m_r}_{-1,xr}}{(m_1+a-r+1)^{k_1} (m_2+a-r+2)^{k_2}\ldots (m_r+a)^{k_r}}.
\end{equation}
\begin{proof}
Using (\ref{HLerch1}) and (\ref{rbStirling}), we have
$$\sum_{n=0}^{\infty}B_{n,a}^{(k_1,k_2,\ldots, k_r)}(x)\frac{t^n}{n!}={\Phi}_{(k_1,k_2,\ldots, k_r)}(1-e^{-t},a)e^{xrt}\qquad\qquad\qquad\qquad\qquad\qquad\qquad\qquad\qquad\qquad\qquad\qquad\qquad\qquad$$
\begin{align*}
&=\sum_{ 0\le m_1\le m_2\le\ldots \le m_r }\frac{m_r!}{(m_1+a-r+1)^{k_1} (m_2+a-r+2)^{k_2}\ldots (m_r+a)^{k_r}}\frac{e^{xrt}(e^{-t}-1)^{m_r}}{(-1)^{m_r}m_r!}\\
&=\sum_{ 0\le m_1\le m_2\le\ldots \le m_r }\frac{m_r!}{(m_1+a-r+1)^{k_1} (m_2+a-r+2)^{k_2}\ldots (m_r+a)^{k_r}}\sum_{n=m_r}^{\infty}\braced{n}{m_r}_{-1,xr}\frac{t^n}{n!}\\
&=\sum_{n=0}^{\infty}\sum_{ 0\le m_1\le m_2\le\ldots \le m_r\le n }\frac{m_r!\braced{n}{m_r}_{-1,xr}}{(m_1+a-r+1)^{k_1} (m_2+a-r+2)^{k_2}\ldots (m_r+a)^{k_r}}\frac{t^n}{n!}
\end{align*}
Comparing the coefficients of $\frac{t^n}{n!}$ completes the proof of the theorem.
\end{proof}
\end{thm}

Note that (\ref{rbStirling}) implies 
$$\braced{n}{m_r}_{-1,0}=(-1)^{n+m_r}\braced{n}{m_r}.$$ 
Hence, as a direct consequence of Theorem \ref{res1}, we have the following corollary.

\begin{cor} 
The Hurwitz-Lerch type multi poly-Bernoulli numbers equal
\begin{equation}\label{cor1}
B_{n,a}^{(k_1,k_2,\ldots, k_r)}=\sum_{0\le m_1\le m_2\le\ldots \le m_r\le n}\frac{(-1)^{n+m_r}m_r!\braced{n}{m_r}}{(m_1+a-r+1)^{k_1} (m_2+a-r+2)^{k_2}\ldots (m_r+a)^{k_r}}.
\end{equation}
\end{cor}

When $r=1$, equation (\ref{cor1}) gives
$$B_{n,a}^{(k_1)}=\sum_{0\le m_1\le n}\frac{(-1)^{n+m_1}m_1!\braced{n}{m_1}}{(m_1+a)^{k_1}}$$
which is exactly the explicit formula for Hurwitz type poly-Bernoulli numbers in Theorem 2.1 of \cite{Cenkci}.

\bigskip
\begin{flushleft}
{\bf Roberto B. Corcino}\\
Cebu Normal University\\
Cebu City, Philippines\\
e-mail: rcorcino@yahoo.com
\end{flushleft}

\begin{flushleft}
{\bf Hassan Jolany}\\
Universit\'e des Sciences et Technologies de Lille\\
UFR de Math\'ematiques\\
Laboratoire Paul Painlev\'e\\
CNRS-UMR 8524 59655 Villeneuve d'Ascq Cedex/France\\
e-mail: hassan.jolany@math.univ-lille1.fr
\end{flushleft}

\begin{flushleft}
{\bf Cristina B. Corcino}\\
Cebu Normal University\\
Cebu City, Philippines\\
e-mail: cristinacorcino@yahoo.com
\end{flushleft}

\begin{flushleft}
{\bf Takao Komatsu}\\
School of Mathematics and Statistics\\
Wuhan University\\ 
Wuhan 430072 China\\
e-mail: komatsu@whu.edu.cn
\end{flushleft}


\begin{thebibliography}{}
\bibitem{Araci} S. Araci, M. Acikgoz and E. Sen, On the extended Kim's $p$-adic $q$-deformed fermionic integrals in the p-adic integer ring, {\it J. Number Theory}, {\bf 133} (2013), 3348--3361. 
\bibitem{Bayad} A. Bayad and Y. Hamahata, Arakawa-Kaneko $L$-functions and generalized poly-Bernoulli polynomials, {\it J. Number Theory}, {\bf 131} (2011), 1020--1036.
\bibitem{Bayad2} A. Bayad and Y. Hamahata, Multiple polylogarithms and multi-poly-Bernoulli polynomials, {\it Funct. Approx. Comment. Math.}, {\bf 46} (2012), 45--61.
\bibitem{Benyi} B. Be\'nyi, Advances in Bijective Combinatorics, {\it Ph.D. Thesis}, 2014.
\bibitem{Brewbaker} C. Brewbaker, A Combinatorial Interpretation of the Poly-Bernoulli Numbers and Two Fermat Analogues, {\it Integers}, {\bf 8} (2008), \#A02.
\bibitem{Candel} B. Candelpergher and M. A. Coppo, A new class of identities involving Cauchy numbers, harmonic numbers and zeta values, {\it Ramanujan J.}, {\bf 27} (2012), 305--328.
\bibitem{Cenkci} M. Cenkci and P. T. Young, Generalizations of Poly-Bernoulli and Poly-Cauchy Numbers, {\it Eur. J. Math.}, Published Online on 01 September 2015, DOI 10.1007/s40879-015-0071-3.
\bibitem{Comtet} L. Comtet, Advanced Combinatorics, {\it D. Reidel Publishing Company}, 1974.
\bibitem{Coppo} M-A. Coppo and B. Candelpergher, The Arakawa-Kaneko Zeta Function, {\it Ramanujan J.}, {\bf 22} (2010), 153--162.
\bibitem{CORCA} R. B. Corcino, C. B. Corcino and R. Aldema, Asymptotic Normality of the $(r,\beta)$-Stirling Numbers, {\it Ars Combin.}, {\bf 81} (2006), 81--96.
\bibitem{CorJol2} R. B. Corcino, H. Jolany, M. Aliabadi and M. R. Darafsheh,  A Note on Multi Poly-Euler Numbers and Bernoulli Polynomials, {\it General Mathematics}, {\bf 20} (2-3) (2012), 122--134 (ROMANIA).
\bibitem{Hamahata} Y. Hamahata, Poly-Euler Polynomials and Arakawa-Kaneko Type Zeta Functions, {\it Funct. Approx. Comment. Math.}, {\bf 51}(1) (2014), 7-22.
\bibitem{Imatomi} K. Imatomi, M. Kaneko and E. Takeda, Multi-Poly-Bernoulli Numbers and Finite Multiple Zeta Values, {\it J. Integer Seq.}, {\bf 17} (2014), Article 14.4.5.
\bibitem{Jang} L. Jang, T. Kim, and H. K. Pak, A note on $q$-Euler and Genocchi numbers, {\it Proc. Japan Acad. Ser. A Math. Sci.}, {\bf 77} (2001), 139--141.
\bibitem{Jolany1} H. Jolany, R. E. Alikelaye and S. S. Mohamad, Some Results on the Generalization of Bernoulli, Euler and Genocchi Polynomials, {\it Acta Univ. Apulensis Math. Inform.}, {\bf 27} (2011), 299--306.
\bibitem{CorJol1} H. Jolany, R. B. Corcino and T. Komatsu, More Properties of Multi Poly-Euler Polynomials, {\it Bol. Soc. Mat. Mex.}, {\bf 21} (2015), 149--162.
\bibitem{Jolany3} H. Jolany, M.R. Darafsheh, R.E. Alikelaye, Generalizations of Poly-Bernoulli Numbers and Polynomials, {\it Int. J. Math. Comb.}, {\bf 2} (2010), 7--14.
\bibitem{Kaneko} M. Kaneko, Poly-Bernoulli numbers, {\it J. Th\'eor. Nombres Bordeaux}, {\bf 9} (1997), 221--228.
\bibitem{Kim} T. Kim, $q$-Generalized Euler numbers and polynomials, Russ. J. Math. Phys. {\bf 13}(3) (2006), 293--298.
\bibitem{Lee} D. W. Lee, On Multiple Appell Polynomials, {\it Proc. Amer. Math. Soc.}, {\bf 139} (2011), 2133--2141. 
\bibitem{Sasaki} Y. Ohno and Y. Sasaki, On the parity of poly-Euler numbers, {\it RIMS Kokyuroku Bessatsu}, {\bf B32} (2012), 271--278.
\bibitem{Shohat} J. Shohat, The Relation of the Classical Orthogonal Polynomials to the Polynomials of Appell, {\it Amer. J. Math.}, {\bf 58} (1936), 453-464
\bibitem{Toscano} L. Toscano, Polinomi Ortogonali o Reciproci di Ortogonali Nella classe di Appell, {\it Le Matematiche}, {\bf 11} (1956), 168--174.
\end{thebibliography}
\end{document}